\documentclass[11pt]{article}
\usepackage{geometry}
\usepackage{dsfont}

\usepackage{psfrag}

\usepackage{calrsfs}
\usepackage{mathrsfs}
\usepackage{pdfsync}
\usepackage{amsmath, amsthm, amssymb, amsfonts}
\usepackage[shortlabels]{enumitem}

\usepackage{url}
\usepackage{float}

\usepackage[usenames]{color}

\usepackage{tikz}

\usetikzlibrary{arrows,decorations.pathmorphing,backgrounds,positioning,fit,automata}

\usetikzlibrary{arrows}
\usetikzlibrary{petri}
\usetikzlibrary{topaths}

\usepackage{indentfirst,calc,euscript}
\usepackage{setspace}
\usepackage[reals]{layout}
\usepackage{xr}
\usepackage{amscd}

\usepackage{sgame}
\usepackage{subfigure}
\usepackage{slashbox}

\usepackage[colorlinks=true,breaklinks=true,bookmarks=true,urlcolor=blue,
     citecolor=blue,linkcolor=blue,bookmarksopen=false,draft=false]{hyperref}

\newcommand{\ignore}[1]{}

\newtheorem{theorem}{Theorem}

\theoremstyle{definition}

\newtheorem{case}{Case}

\newtheorem{property}{Property}

\numberwithin{equation}{section}

\numberwithin{theorem}{section}

\newcommand{\m}{\mathbb}
\newcommand{\sgn}{\operatorname{Sgn}}

\author{Bruno Ziliotto}
\title{Stochastic homogenization of nonconvex Hamilton-Jacobi equations: a counterexample}
\begin{document}
\maketitle
\begin{abstract}
This paper provides a counterexample to Hamilton-Jacobi homogenization in the nonconvex case, for general stationary ergodic environments. 
\end{abstract}


\bibliographystyle{plain}  





\section*{Introduction}
We consider Hamilton-Jacobi equations of the form 
\begin{equation} \label{HJgen}
\left\{ \begin{array}{l}
\partial_t u(x,t,\omega) + H(Du(x,t,\omega), \frac{x}{\epsilon},\omega)=0 \qquad {\rm in}\; \m{R}^n\times (0,+\infty)\\
u(x,0)= 0 \qquad {\rm in}\; \m{R}^n
\end{array} \right.
\end{equation}
where $n \geq 1$, $\epsilon>0$, the Hamiltonian $H(p,x,\omega)$ is coercive and Lipschitz in $p$, and depends on a random element $\omega$ lying in a probability space $(\Omega,\mathcal{F},\m{P})$. The law of $\omega \rightarrow H(.,\omega)$ is assumed to be stationary and ergodic. Under these assumptions, it is well-known that the above system admits a unique viscosity solution $u^{\epsilon}$, which is measurable with respect to $\omega$. A central question in the literature is to study the convergence properties of $u^{\epsilon}$, as $\epsilon$ goes to 0. The case where the Hamiltonian is periodic in space has been studied by Lions, Papanicolaou and Varadhan \cite{LPV86} (see also Evans \cite{E92}).
The first result in the stochastic case has been obtained by Souganidis \cite{S99} and Rezakhanlou and Tarver \cite{RT00}. They have proved independently that when $H$ is convex with respect to $p$, then $u^{\epsilon}$ converges $\m{P}$-almost surely to the unique solution of a system of the form
$$
\left\{ \begin{array}{l}
\partial_t u(x,t) + \bar{H}(Du(x,t))=0 \qquad {\rm in}\; \m{R}^n\times (0,+\infty)\\
u(x,0)= 0 \qquad {\rm in}\; \m{R}^n
\end{array} \right.
$$
where $\bar{H}$ is the \textit{effective Hamiltonian}.
This result has been extended to various frameworks, still under the assumption that the Hamiltonian is convex in $p$ (see \cite{LS05,KRV06,LS10,LS03,S09,AS131,AT141}). Quantitative results about the speed of convergence have been obtained in \cite{ACS14,MN12,AC15a}.
 The question of the homogenization of Hamilton-Jacobi equations in the general case where $H$ is not convex in $p$ had remained open until now, and is regularly mentioned in the literature (see for instance \cite{LS05,LS10,K07,AC15b,ATY15,ATY152}). A few particular cases have been treated, for example the case of level-set convex Hamiltonians (see Armstrong and Souganidis \cite{AS132}), the case where the law of $H$ is invariant by rotation (this is a direct consequence of Fehrman \cite[Theorem 1.1]{F14}), the 1-dimensional case (see Armstrong, Tran and Yu \cite{ATY152} and Gao \cite{G15}), and the case where the law of $H$ satisfies a finite range condition (see Armstrong and Cardaliaguet \cite{AC15b}).
 \\
 
In this paper, we give a negative answer to this question. Indeed, we provide an example of a Hamilton-Jacobi equation of the form (\ref{HJgen}) in the two-dimensional case ($n=2$), such that $\m{P}$-almost surely, $u_{\epsilon}(0,1,\omega)$ does not converge when $\epsilon$ goes to 0: there is no stochastic homogenization for this equation. In this example, the Hamiltonian $H$ satisfies all the standard assumptions of the literature, except the convexity with respect to $p$. The Hamilton-Jacobi equation is associated to a zero-sum differential game. A formal description of such a game is unnecessary to prove the main theorem. Instead, an informal discussion is provided.

 The paper is organized as follows. Section 1 is dedicated to the construction of the example, and Section 2 proves the main result of the paper. Section 3 examines how the Hamiltonian of the example correlates distant regions of space.

\section{Construction of the example}
\subsection{The weight function}
For all $n \geq 1$, equip $\m{R}^n$ with the Euclidean norm $\left|.\right|$. Let $E$ be the set of 1-Lipschitz mappings from $\m{R}^2$ to $[1,2]$. Let us build a probability measure on $E$ in the following way.

Let $(T_k)$ be the sequence defined for $k \geq 1$ by $T_k= 4^{k}$. 
\\
Let $(X^j_{k,l,m})_{(j,k,l,m) \in \left\{1,2\right\} \times \m{N}^*\times \m{Z}^2}$ be a sequence of independent random variables defined on a probability space $(\Omega,\mathcal{F},\m{P})$, such that for all $(j,k,l,m) \in \left\{1,2\right\} \times \m{N}^* \times \m{Z}^2$, $X^j_{k,l,m}$ follows a Bernoulli of parameter $T_k^{-2}$. 

Let $\omega \in \Omega$. The mapping $c_{\omega}:\m{R}^2 \rightarrow [1,2] \in E$ is built in three phases.

\vspace{0.5cm}
\textbf{Phase 1}
\vspace{0.5cm}
\\
The mapping $c^1_{\omega}:\m{R}^2 \rightarrow [0,2]$ is built through the following step-by-step procedure:
\\
\begin{itemize}
\item
Step $k=0$: take $c^1_{\omega}:=0$ as the initial distribution of weights. 
\item
Step $k \geq 1$: for each $(l,m) \in \m{Z}^2$ such that $X^{1}_{k,l,m}(\omega)=1$, consider the horizontal segment centered on $(l,m)$, with length $10T_k$, which shall be called \textit{green segment of length $10 T_k$}. For each $(x_1,x_2) \in \m{R}^2$ that lies in the segment, set $c^1_{\omega}(x_1,x_2):=1$ (note that $c^1_{\omega}(x_1,x_2)$ may have already been defined as being 1 at an earlier step).
\end{itemize}
At the end of Phase 1, we have a map $c^1_{\omega}: \m{R}^2 \rightarrow [0,2]$. Then, go to Phase 2:

\vspace{0.5cm}
\textbf{Phase 2}
\vspace{0.5cm}
\\
The mapping $c^2_{\omega}:\m{R}^2 \rightarrow [0,2]$ is built through the following step-by-step procedure:
\begin{itemize}
\item
Step $k=0$: take $c^2_{\omega}:=c^1_{\omega}$ as the initial distribution of weights.
\item
Step $k \geq 1$: for each $(l,m) \in \m{Z}^2$ such that $X^{2}_{k,l,m}(\omega)=1$, consider the vertical segment centered on $(l,m)$, with length $10T_k$, which shall be called \textit{red segment of length $10 T_k$}. For each $(x_1,x_2)$ that lies in the segment, proceed as follows:
\begin{itemize}
\item
If a green segment of size $10 T_{k'}$ with $k' \geq k$ is at a distance strictly smaller than 1 of $(x_1,x_2)$, $c^{2}_{\omega}(x_1,x_2)$ is not modified. 
\item
Otherwise, set $c^{2}_{\omega}(x_1,x_2):=2$. Moreover, for all $x \in (x_1-1,x_1+1)\setminus\left\{x_1\right\} \times \left\{x_2\right\}$, set $c^2_{\omega}(x):=0$.  
\end{itemize}
A key feature of Phase 2 is that, whenever a red segment intersects a green segment, the intersection ``turns green" if the red segment's length is smaller than the green segment's one, and ``turns red" otherwise. 
\end{itemize}

A green segment is \textit{complete} if all its elements $(x_1,x_2)$ satisfy $c^2_{\omega}(x_1,x_2)=1$: in other words, it is not intersected by a larger red segment. In the same vein, a red segment is \textit{complete} if all its elements $(x_1,x_2)$ satisfy $c^2_{\omega}(x_1,x_2)=2$.

\vspace{0.5cm}
\textbf{Phase 3}
\vspace{0.5cm}
\\
Define $c_{\omega}:\m{R}^2 \rightarrow [1,2]$ by
$$
c_\omega(x):= \max \left(\sup_{y\in \m{R}^2, \; c^2_\omega(y)>0} \left\{ c^2_\omega(y)-\left|x-y\right| \right\},1 \right).
$$
Note that if $c^2_{\omega}(x) \in \left\{1,2\right\}$, then $c_{\omega}(x)=c^2_{\omega}(x)$. Moreover, for all $\omega \in \Omega$, $c_{\omega}$ is 1-Lipschitz.
By construction, the law of the random variable $\omega \rightarrow c_{\omega}$ is $\m{Z}^2$-invariant. Let us check that it is ergodic, that is, for all event $A$ which is $\m{Z}^2$-invariant, $\m{P}(A)=0$ or $\m{P}(A)=1$. Let $A$ be a $\m{Z}^2$-invariant event. 

For $n \in \m{N}^*$, denote by $\mathcal{F}_n$ the $\sigma$-algebra generated by the random variables 
\\
${(X^j_{k,l,m})}_{(j,k,l,m) \in \left\{1,2\right\} \times \m{N}^* \times [-n,n]^2}$. Let $\epsilon>0$. 
There exists $n \geq 1$ and $A_n \subset A$ such that $A_n$ is $\mathcal{F}_n$-measurable and $\m{P}(A_n) \geq \m{P}(A)-\epsilon$. Let $A'_n$ be the translation of $A_n$ with respect to the vector $(0,2n+1)$. Because $A$ is $\m{Z}^2$-invariant, $A'_n$ is a subset of $A$. Because the law of $\omega \rightarrow c_{\omega}$ is $\m{Z}^2$-invariant, the equality $\m{P}(A'_n)=\m{P}(A_n)$ holds. By construction, the events $A_n$ and $A'_n$ are independent, thus
\begin{eqnarray*}
\m{P}(A \cap A) &\leq& \m{P}(A_n \cap A'_n)+2 \epsilon
\\
&=&\m{P}(A_n) \m{P}(A'_n)+2\epsilon
\\
&\leq& \m{P}(A)^2+5\epsilon.
\end{eqnarray*}
It follows that $\m{P}(A) \leq \m{P}(A)^2$, thus $\m{P}(A)=0$ or $\m{P}(A)=1$: the law of $\omega \rightarrow c_{\omega}$ is ergodic.
\\

\subsection{Main result}
Let $l: \m{R}^2 \times [-1,1]^2 \times \Omega \rightarrow [1,2]$ defined by 
\begin{equation*}
\forall (x,a,\omega) \in \m{R}^2 \times [-1,1]^2 \times \Omega, \quad l(x,a,\omega):=c_{\omega}(x)+10 |a_1|,
\end{equation*}
and $H: \m{R}^2 \times \m{R}^2 \times \Omega$ defined by
$$
\forall (p,x,\omega) \in \m{R}^2 \times \m{R}^2 \times \Omega, \quad H(p,x, \omega):= \max_{a\in [-1,1]^2}\min_{b\in [-1,1]^2} \{-l(x,a,\omega)-p\cdot (2a+b)\}.
$$
Note that for all $\omega \in \Omega$, $(x,a) \rightarrow l(x,a,\omega)$ is 10-Lipschitz, and $H$ is coercive in $p$, uniformly in $x$ and $\omega$. 
For $\epsilon>0$, consider the following Hamilton-Jacobi equation:
\begin{equation} \label{HJ1}
\left\{ \begin{array}{l}
\partial_t u(x,t,\omega) + H(Du(x,t,\omega), \frac{x}{\epsilon},\omega)=0 \qquad {\rm in}\; \m{R}^2\times (0,+\infty)\\
u(x,0)= 0 \qquad {\rm in}\; \m{R}^2,
\end{array} \right.
\end{equation}
where $\partial_t u$ and $Du$ are, respectively, the temporal derivative and the gradient of $u$.
\begin{theorem} \label{main}
Let $u^{\epsilon}$ be the solution of (\ref{HJ1}). Then
\begin{equation*}
\liminf_{\epsilon \rightarrow 0} u^{\epsilon}(0,1,\omega)=1 \quad \text{and} \quad \limsup_{\epsilon \rightarrow 0} u^{\epsilon}(0,1,\omega)=2 \quad \m{P}-\text{almost surely}.
\end{equation*}
Consequently, there is no stochastic homogenization for the above Hamilton-Jacobi equation.
\end{theorem}
\section{Proof of the theorem}
\subsection{A change of variables}
For the proof, it is easier to consider the following system:
\begin{equation} \label{HJB}
\left\{ \begin{array}{l}
\partial_t u(x,t,\omega) + H(Du(x,t,\omega), x,\omega)=0 \qquad {\rm in}\; \m{R}^2\times (0,+\infty)\\
u(x,0)= 0 \qquad {\rm in}\; \m{R}^2.
\end{array} \right.
\end{equation}
The solution $u$ of the above system satisfies the relation
\begin{equation} \label{scaling}
\forall \epsilon >0 \quad \forall t \in \m{R}_+ \quad u^{\epsilon}(0,t,\omega)=\epsilon u(0,t/\epsilon,\omega) \quad \m{P}-\text{almost surely}.
\end{equation}
  In the remainder of the paper, we prove that
\begin{equation*}
\liminf_{T \rightarrow+\infty} \frac{1}{T}u(0,T,\omega)=1 \quad \text{and} \quad \limsup_{T \rightarrow+\infty} \frac{1}{T} u(0,T,\omega)=2 \quad \m{P}-\text{almost surely}.
\end{equation*}
By (\ref{scaling}), this implies Theorem \ref{main}.

\subsection{Intuition of the result}
The Hamilton-Jacobi equation of the previous subsection can be associated to a zero-sum stochastic differential game. To prove Theorem \ref{main}, giving a formal description of this game is unnecessary. Nonetheless, in order to have a better understanding of the example, an informal description is given, in which technical details are avoided, notably concerning the definition of strategies. 

Let $\omega \in \Omega$ and $T>0$. The game starts at the origin $(0,0)$, and has a duration $T$. Player 1 (resp. 2) aims at minimizing (resp. maximizing) the total cost between time 0 and time $T$, given by $\int_{0}^T l(x(t),a(t),\omega) dt$, where $a(t)$ is the control used by Player 1 at time $t$. The cost function is the sum of the weight function $c_{\omega}$, which is space-dependent, and the function $a\rightarrow 10 \left| a_1 \right|$, which heavily penalizes the controls of Player 1 that have a nonzero horizontal component. 

The dynamics of the state is such that if Player 1 chooses a control $a \in [-1,1]^2$ and Player 2 chooses a control $b \in [-1,1]^2$, then the state moves in the direction $2a+b$. Thus, Player 1 can control the state and bring it wherever he wants in linear time. The value of the game with duration $T$ coincides with $u(0,T,\omega)$, where $u$ is the solution of the system (\ref{HJB}). 

Ideally, Player 1 would like to bring the state in a region where the weight function takes small values, and force the state to stay there. The issue is that forcing the state to stay in this region may turn out to be very costly for him, if this requires the use of horizontal controls. For all $x \in \m{R}$, denote by $\lfloor x \rfloor$ the integer part of $x$. The construction of $(c_{\omega})_{\omega \in \Omega}$ has been made such that for all $\epsilon>0$, there exist two positive probability events $\Omega_1$ and $\Omega_2$ such that the following properties hold:
\begin{property}
For all $\omega \in \Omega_1$, there exists a sequence $(n_k(\omega))$ going to infinity such that for all $k \geq 1$, there exists a complete green segment of length $10 T_{n_k(\omega)}$ whose center is at a distance smaller or equal to $\lfloor \epsilon T_{n_k(\omega)} \rfloor$ from the origin. 
\end{property}
\begin{property}
For all $\omega \in \Omega_2$, there exists a sequence $(n'_k(\omega))$ going to infinity such that for all $k \geq 1$, there exists a complete red segment of length $10 T_{n_k(\omega)}$ whose center is at a distance smaller or equal to $\lfloor \epsilon T_{n'_k(\omega)} \rfloor$ from the origin.
\end{property}
Let $k \geq 1$ and $\omega \in \Omega_1$. Consider the game with duration $T_{n_k(\omega)}$. Player 1 can force the state to go the center of the complete green segment, within a length of time smaller or equal to $\lfloor \epsilon T_{n_k(\omega)} \rfloor$. Then, he can force it to stay in the segment until the end of the game, by making use only of vertical controls, which are costless. Thus, for $\epsilon$ small enough, the normalized value $u(0,T_{n_k(\omega)},\omega)/T_{n_k(\omega)}$ of the game with duration $T_{n_k(\omega)}$ is close to 1. 
\\

Let $k \geq 1$ and $\omega \in \Omega_2$. Now consider the game with duration $T_{n'_k(\omega)}$. Player 2 can adopt the following strategy: always play in the direction of the complete red segment. If Player 1 does not use horizontal controls, after a length of time smaller or equal to $\lfloor \epsilon T_{n'_k(\omega)} \rfloor$, the state will be in the complete red segment, where the weight function takes the value 2. If Player 1 uses horizontal controls to counter Player 2's controls, he faces a huge penalty cost. It follows that for $\epsilon$ small enough, the normalized value $u(0,T_{n'_k(\omega)},\omega)/T_{n'_k(\omega)}$ of the game with duration $T_{n'_k(\omega)}$ is close to 2. 
\\

Because the law of $H$ is ergodic, these two arguments prove the theorem. 
The next two subsections are dedicated to the formal proof. Subsection \ref{super} exploits Property 1 to build supersolutions of (\ref{HJB}). Subsection \ref{sub} exploits Property 2 to build subsolutions of (\ref{HJB}).

\subsection{Supersolutions of the Hamilton-Jacobi equation} \label{super}
Let $\epsilon \in (0,1/20]$ and $k \geq 1$. Consider the event $B_k$ ``the center of a complete green segment of length $10T_k$ is at a distance smaller or equal to $\lfloor \epsilon T_k \rfloor$ from the origin''. A sufficient condition for $B_k$ to be realized is that the two following events $C_k$ and  $D_k$ are realized:
\begin{itemize}
\item
At step $k$ of Phase 1, a point at a distance smaller or equal to $\lfloor \epsilon T_k \rfloor$ from the origin has been selected by the Bernoulli random variable. 
\item
The green segment of length $10T_k$ centered on this point is complete, that is, it is not intersected by a red segment of length strictly larger than $10T_k$.
\end{itemize}
We have 
\begin{equation*}
\m{P}(C_k) \geq 1-(1-T_k^{-2})^{{(\lfloor \epsilon T_k \rfloor})^2}. 
\end{equation*}
For $k' \geq k+1$, the probability that no red segment of length $10T_{k'}$ intersects the green segment is greater than $(1-T_{k'}^{-2})^{ (10T_{k'}+1) (10T_k+1)}$. Thus,
\begin{equation*}
\m{P}(D_k|C_k) \geq \prod_{k' \geq k+1}(1-T_{k'}^{-2})^{(10T_{k'}+1) (10T_k+1)}.
\end{equation*}
Consequently,
\begin{eqnarray*}
\m{P}(B_k)&\geq&\m{P}(D_k|C_k) \m{P}(C_k) 
\\
&\geq& \left[1-(1-T_k^{-2})^{{(\lfloor \epsilon T_k \rfloor})^2}\right] \left[\prod_{k' \geq k+1}(1-T_{k'}^{-2})^{(10T_{k'}+1) (10T_k+1)} \right].
\end{eqnarray*}
We have $\displaystyle \lim_{k \rightarrow+\infty}\left[1-(1-T_k^{-2})^{{(\lfloor \epsilon T_k \rfloor})^2}\right]>0$. Moreover,
\begin{equation*}
\sum_{k' \geq k+1} T_{k'}^{-2} T_{k'} T_k =T_k \sum_{k' \geq k+1} T_{k'}^{-1},
\end{equation*}
thus $\liminf_{k \rightarrow+\infty} \left[\prod_{k' \geq k+1}(1-T_{k'}^{-2})^{(10T_{k'}+1) (10T_k+1)} \right] >0$. We deduce that 
\\
$\liminf_{k \rightarrow +\infty} \m{P}(B_k)>0$. Thus, there exists a positive probability event $\Omega_1 \subset \Omega$ such that for all $\omega \in \Omega_1$, the events $(B_k)_{k \geq 1}$ occur infinitely often. 
\\
 
Let $\omega \in \Omega_1$ and $k \geq 1$ such that $B_k$ is realized. Let $(X_1(\omega),X_2(\omega))$ be the coordinates of the center of the associated complete green segment of length $10 T_{k}$. In what follows, for simplicity, we omit the dependence in $\omega$. 
\\

Define $u_{+}: \m{R}^2 \times (0,T_k) \rightarrow \m{R}$ by 
\begin{equation*}
u_+(x,t):= 3\left|x_2-X_2\right|+t+(\left|x_1-X_1\right|-5T_k+2t)_+,
\end{equation*}
where for all real-valued function $f$, $(f)_+:=\max(f,0)$. 

Let us prove that $u_+$ is a supersolution of the system (\ref{HJB}). Let $(x,t) \in \m{R}^2 \times (0,T_k)$. We distinguish between the following cases:
\begin{case} $x_2 \neq X_2$ 
\end{case}
Let $\phi$ be a smooth function such that $\phi(x,t)=u_+(x,t)$ and $\phi \leq u_+$ on a neighborhood of $(x,t)$.
 Then $\partial_t \phi(x,t) \geq 1$ and $\partial_{x_2} \phi(x,t)=3\sgn(x_2-X_2)$, where $\sgn$ is the sign function. Moreover, $\left|\partial_{x_1} \phi(x,t)\right| \leq 1$.
For $a=(0,-\sgn(x_2-X_2))$, for all $b \in [-1,1]^2$, 
\begin{eqnarray*}
-l(x,a) - D \phi(x,t) \cdot (2a+b) \geq -2+3(2-1)-1=0,
\end{eqnarray*}
thus
\begin{eqnarray*}
\partial_t \phi(x,t)+\max_{a \in [-1,1]^2} \min_{b \in [-1,1]^2} \left\{-l(x,a) - D\phi(x,t) \cdot (2a+b) \right\} &\geq& 0.
\end{eqnarray*}
\begin{case} $\left|x_1-X_1\right|-5T_k+2t \leq 0$ and $x_2 = X_2$
\end{case}

Let $\phi$ be a smooth function such that $\phi(x,t)=u_+(x,t)$ and $\phi \leq u_+$ on a neighborhood of $(x,t)$. Then $\partial_t \phi(t,x)-\left|\partial_{x_1}\phi(t,x)\right| \geq 1$. Let $b \in [-1,1]^2$ and $a:=(0,-b_2/2)$. Because $\left|x_1-X_1\right|-5T_k+2t \leq 0$, we have $l(x,a)=1$. Then
\begin{equation*}
-l(x,a) - D\phi(x,t) \cdot (2a+b) \geq -1-\left|\partial_{x_1}\phi(t,x)\right|,
\end{equation*}
and
\begin{eqnarray*}
\partial_t \phi(x,t)+\max_{a \in [-1,1]^2} \min_{b \in [-1,1]^2} \left\{-l(x,a) -  D\phi(x,t) \cdot (2a+b) \right\} &\geq& -1+1=0.
\end{eqnarray*}
\begin{case}{$\left|x_1-X_1\right|-5T_k+2t>0$ and $x_2=X_2$}
\end{case}
Let $\phi$ be a smooth function such that $\phi(x,t)=u_+(x,t)$ and $\phi \leq u_+$ on a neighborhood of $(x,t)$. Then $\partial_t \phi(x,t)=3$ and $\left|\partial_{x_1}\phi(x,t)\right|\leq 1$. Let $b \in [-1,1]^2$, and $a:=(0,-b_2/2)$. Then
\begin{eqnarray*}
-l(x,a) - D \phi(x,t) \cdot (2a+b) \geq -2-1=-3.
\end{eqnarray*}
We deduce that
\begin{eqnarray*}
\partial_t \phi(x,t)+\max_{a \in [-1,1]^2} \min_{b \in [-1,1]^2} \left\{-l(x,a) - D \phi(x,t) \cdot (2a+b) \right\} &\geq& 0.
\end{eqnarray*}
Consequently, $u_+$ is a supersolution of the system (\ref{HJB}). 
Comparison principle (see Crandall, Ishii and Lions \cite{CIL92}) implies that for all $\omega \in \Omega_1$ 
\begin{equation*}
u(0,T_k,\omega) \leq u_+(0,T_k)=3\left|X_2\right|+T_k \leq 3 \lfloor \epsilon T_k \rfloor +T_k.
\end{equation*}
We deduce that for all $\omega \in \Omega_1$, 
\begin{equation*}
\liminf_{T \rightarrow+\infty} T^{-1} u(0,T,\omega)=1.
\end{equation*}
The map $u$ is uniformly Lipschitz with respect to $x$, and the law of $H$ is ergodic. A well-known consequence is that the random variable $\liminf_{T \rightarrow+\infty} T^{-1} u(0,T,\omega)$ is $\m{P}$-almost surely constant. This implies that $\m{P}$-almost surely,
\begin{equation*}
\liminf_{T \rightarrow+\infty} T^{-1} u(0,T,\omega)=1.
\end{equation*}
\subsection{Subsolutions of the Hamilton-Jacobi equation} \label{sub}
Let $\epsilon \in (0,1/20]$ and $k \geq 1$. Consider the event $B'_k$ ``the center of a complete red segment of length $10T_k$ is at a distance smaller or equal to $\lfloor \epsilon T_k \rfloor$ from the origin''. A sufficient condition for $B'_k$ to be realized is that the two following events $C'_k$ and  $D'_k$ are realized:
\begin{itemize}
\item
At step $k$ of Phase 2, a point at a distance smaller or equal to $\lfloor \epsilon T_k \rfloor$ from the origin has been selected by the Bernoulli random variable. 
\item
The red segment of length $10T_k$ centered on this point is complete, that is, it is not intersected by a green segment of length larger or equal to $10T_k$.
\end{itemize}
We have 
\begin{equation*}
\m{P}(C'_k) \geq 1-(1-T_k^{-2})^{{(\lfloor \epsilon T_k \rfloor})^2}. 
\end{equation*}
For $k' \geq k$, the probability that no green segment of length $10T_{k'}$ intersects the red segment is greater than $(1-T_{k'}^{-2})^{(10T_{k'}+1) (10T_k+1)}$. Thus,
\begin{equation*}
\m{P}(D'_k|C'_k) \geq \prod_{k' \geq k}(1-T_{k'}^{-2})^{(10T_{k'}+1) (10T_k+1)}.
\end{equation*}
Consequently,
\begin{eqnarray*}
\m{P}(B'_k)&\geq&\m{P}(D'_k|C'_k) \m{P}(C'_k) 
\\
&\geq& \left[1-(1-T_k^{-2})^{{(\lfloor \epsilon T_k \rfloor})^2}\right] \left[\prod_{k' \geq k}(1-T_{k'}^{-2})^{(10T_{k'}+1) (10T_k+1)} \right].
\end{eqnarray*}

Similar computations as in Subsection \ref{super} show that $\liminf_{k \rightarrow +\infty} \m{P}(B'_k)>0$. Thus, there exists a positive probability event $\Omega_2 \subset \Omega$ such that for all $\omega \in \Omega_2$, the events $(B'_k)_{k \geq 1}$ occur infinitely often. 
\\
 
Let $\omega \in \Omega_2$ and $k \geq 1$ such that $B'_k$ is realized. Let $(X_1(\omega),X_2(\omega))$ be the coordinates of the center of the associated complete red segment of length $10 T_{k}$. In what follows, for simplicity, we omit the dependence in $\omega$. 
\\

Define $u_{-}: \m{R}^2 \times (0,T_k) \rightarrow \m{R}$ by 
\begin{equation*}
u_-(x,t)=2t-3\left|x_1-X_1\right|+(5T_k-\left|x_2-X_2\right|-10t)_-,
\end{equation*}
where for all real-valued function $f$, $(f)_-:=\min(f,0)$. 
Let us prove that $u_-$ is a subsolution of the system (\ref{HJB}). Let $(x,t) \in \m{R}^2 \times (0,T_k)$. We distinguish between the following cases:
\setcounter{case}{0}
\begin{case} $x_1 \neq X_1$ 
\end{case}
Let $\phi$ be a smooth function such that $\phi(x,t)=u_-(x,t)$ and $\phi \geq u_-$ on a neighborhood of $(x,t)$.
 Then $\partial_t \phi(x,t) \leq 2$ and $\partial_{x_1} \phi(x,t)=3\sgn(X_1-x_1)$. Moreover, $\left|\partial_{x_2} \phi(x,t)\right| \leq 1$.
For $b=(\sgn(X_1-x_1),0)$, for all $a \in [-1,1]^2$, 
\begin{eqnarray*}
-l(x,a) - D \phi(x,t) \cdot (2a+b) \leq -1-10|a_1|+3(2|a_1|-1)+2 \leq -2,
\end{eqnarray*}
thus
\begin{eqnarray*}
\partial_t \phi(x,t)+\max_{a \in [-1,1]^2} \min_{b \in [-1,1]^2} \left\{-l(x,a) - D\phi(x,t) \cdot (2a+b) \right\} &\leq& 0.
\end{eqnarray*}

\begin{case}{$5T_k-\left|x_2-X_2\right|-10t \geq 0$ and $x_1=X_1$}
\end{case}
The key point is that in this case, $c_{\omega}(x)=2$.
Let $\phi$ be a smooth function such that $\phi(x,t)=u_-(x,t)$ and $\phi \geq u_-$ on a neighborhood of $(x,t)$. Let $a \in [-1,1]^2$. We have $\partial_t \phi(x,t)+2\left| \partial_{x_2}\phi(x,t)\right| \leq 2$ and $\left| \partial_{x_1}\phi(x,t)\right| \leq 3$. Let $b:=0$. Then
\begin{eqnarray*}
-l(x,a) - D \phi(x,t) \cdot (2a+b) &\leq& -2-10|a_1|-2\partial_{x_1}\phi(x,t)a_1-2\partial_{x_2} \phi(x,t)a_2 
\\
&\leq& -2+2\left|\partial_{x_2} \phi(x,t)\right|.
\end{eqnarray*}
We deduce that
\begin{eqnarray*}
\partial_t \phi(x,t)+\max_{a \in [-1,1]^2} \min_{b \in [-1,1]^2} \left\{-l(x,a) - D \phi(x,t) \cdot (2a+b) \right\} &\leq& 0.
\end{eqnarray*}
\begin{case}{$5T_k-\left|x_2-X_2\right|-10t < 0$ and $x_1=X_1$}
\end{case}
Let $\phi$ be a smooth function such that $\phi(x,t)=u_-(x,t)$ and $\phi \geq u_-$ on a neighborhood of $(x,t)$. Then $\partial_t \phi(x,t)=-8$, and $\left|\partial_{x_1}\phi(x,t)\right|+ \left|\partial_{x_2}\phi(x,t)\right|\leq 4$. Thus, for $b=0$ and all $a \in [-1,1]^2$, we have
\begin{eqnarray*}
-l(x,a) - D \phi(x,t) \cdot (2a+b) \leq -1+8=7,
\end{eqnarray*}
thus
\begin{eqnarray*}
\partial_t \phi(x,t)+\max_{a \in [-1,1]^2} \min_{b \in [-1,1]^2} \left\{-l(x,a) - D \phi(x,t) \cdot (2a+b) \right\} &\leq& 0.
\end{eqnarray*}
Consequently, $u_-$ is a subsolution of the system (\ref{HJB}). 
Comparison principle implies that for all $\omega \in \Omega_2$, 
\begin{equation*}
u(0,T_k,\omega) \geq u_-(0,T_k)=2T_k-3\left|X_1\right| \geq 2T_k - 3\lfloor \epsilon T_k \rfloor.
\end{equation*}
We deduce that for all $\omega \in \Omega_2$, 
\begin{equation*}
\limsup_{T \rightarrow+\infty} T^{-1} u(0,T,\omega) = 2.
\end{equation*}
The map $u$ is uniformly Lipschitz with respect to $x$, and the law of $H$ is ergodic. A well-known consequence is that the random variable $\limsup_{T \rightarrow+\infty} T^{-1} u(0,T,\omega)$ is $\m{P}$-almost surely constant. This implies that $\m{P}$-almost surely,
\begin{equation*}
\limsup_{T \rightarrow+\infty} T^{-1} u(0,T,\omega) = 2.
\end{equation*}
\section{Correlation between distant regions of space}
Let us point out that under the law of the Hamiltonian $H$, the correlation between distant regions of space is nonzero. 
In particular, the law of $H$ does not satisfy the finite range condition imposed in \cite{AC15b}. It is natural to ask whether the correlation between two regions of space vanishes as the distance between these two regions goes to infinity. In the literature, several definitions of correlation are considered.
\subsection{A first criterion}
Let $r>0$ and $d>0$. Denote by $\mathcal{O}(r,d)$ the set of pairs of open subsets $(U,V)$ of $\m{R}^2$ such that $\inf_{(x,y) \in U \times V} |x-y| \geq r$, $\sup_{(x,x') \in U^2} |x-x'| \leq d$ and $\sup_{(y,y') \in V^2} |y-y'| \leq d$. Define
\begin{align*}
\rho_1(r,d) & :=\sup \bigg \lbrace \m{P}(E \cap F)-\m{P}(E)\m{P}(F), (U,V) \in \mathcal{O}(r,d),	 \\ \ E \in \mathcal{F}(U), \ F \in \mathcal{F}(V)
 &   \bigg \rbrace,
\end{align*}

where $\mathcal{F}(U)$ is the $\sigma$-algebra generated by the random variables 
\\
$(H(p,x,.))_{(p,x) \in \m{R}^2 \times U}$. 
\\
Let $\gamma>0$. The law of $H$ is polynomially mixing of order $\gamma$ if for all $d>0$, we have $\displaystyle \limsup_{r \rightarrow +\infty} \ r^{\gamma} \rho_1(r,d)<+\infty$ (see Bramson, Zeitouni and Zerner \cite{BZZ06}). We claim that in the example, the law of $H$ is polynomially mixing of order 1. Indeed, fix $d>0$ and let $r>0$. 
Let $(U,V) \in \mathcal{O}(r,d)$ and $(E,F) \in \mathcal{F}(U) \times \mathcal{F}(V)$. Let $A(r,U,V)$ be the event ``there exists neither a green segment nor a red segment of length greater than $r/4$ that crosses either $U$ or $V$". Conditional to $A(r,U,V)$, the events $E$ and $F$ are independent. Moreover, similar computations as in Subsections \ref{super} and \ref{sub} show that 
\begin{equation*}
\limsup_{r \rightarrow +\infty} \ r \sup_{(U,V) \in \mathcal{O}(r,d)} \m{P}(A(r,U,V))<+\infty.
\end{equation*}
These two facts prove that 
$\displaystyle \limsup_{r \rightarrow +\infty} \ r \rho_1(r,d)<+\infty$.

\subsection{A second criterion}
In some papers, like in Yurinskii \cite{Y89}, the correlation at a distance $r>0$ is measured by the quantity $\rho_2(r)$ defined by
\begin{align*}
\rho_2(r) & :=\sup \bigg \lbrace \m{P}(E \cap F)-\m{P}(E)\m{P}(F), \ E \in \mathcal{F}(U), \ F \in \mathcal{F}(V),	 \\
 & \hspace{5mm}  \ U,V \ \text{open subsets of} \ \m{R}^2, \inf_{(x,y) \in U \times V} |x-y| \geq r \bigg \rbrace.
\end{align*}
\\
Note that for all $d>0$, $\rho_2(r) \geq \rho_1(r,d)$. 
Let us prove that $\rho_2(r)$ does not vanish when $r$ goes to infinity. 
Indeed, let $k \geq 2$, $r:=3T_k$ and $x_1>0$. Let $U:=(0,x_1) \times (r,+\infty)$ and $V:=(0,x_1) \times (-\infty,r/2)$. Let $E(x_1)$ be the event ``there exists $a_1$ in $(0,x_1)$ such that there exists a red segment which goes through $(a_1,r)$ and $(a_1,2r)$, and in addition $c_{\omega}((a_1,3r),\omega)<2$''. Let $F$ be the event ``there exists $a_1$ in $(0,x_1)$ such that there exists a red segment which goes through $(a_1,0)$ and $(a_1,r/2)$''. 

Similar computations as in Subsection \ref{super} show that
\begin{equation*}
\lim_{x_1 \rightarrow +\infty} \m{P}(E(x_1))=1.
\end{equation*} 
Let $x_1>0$ such that $\m{P}(E(x_1)) \in [1/2,2/3]$. Assume that $E(x_1)$ is realized. The red segment which goes through $(a_1,r)$ and $(a_1,2r)$ has a length greater than $3T_k > 10 T_{k-1}$. Consequently, it has a length greater than $10 T_{k}$. This implies that the red segment also goes through $(a_1,0)$ and $(a_1,r/2)$. Consequently, $F(x_1)$ is realized. This implies that $\rho_2(r) \geq 1/2-(2/3)^2 = 1/18$. As $k$ has been taken arbitrarily, and $\lim_{k \rightarrow+\infty} T_k=+\infty$, we deduce that $\rho_2(r)$ does not vanish when $r$ goes to infinity. 

A natural question is to ask the following: assuming that $\displaystyle \lim_{r \rightarrow+\infty} \rho_2(r)=0$, is it possible to prove stochastic homogenization?
\section*{Acknowledgments}
I am very grateful to Pierre Cardaliaguet and Sylvain Sorin for their advices and for their careful rereading. I also thank an anonymous referee for his helpful comments.

\bibliography{meancurv2}

\newcommand{\noop}[1]{}
\begin{thebibliography}{10}

\bibitem{AC15b}
S.~Armstrong and P.~Cardaliaguet.
\newblock Stochastic homogenization of quasilinear hamilton-jacobi equations
  and geometric motions.
\newblock {\em preprint arXiv:1504.02045}, 2015.

\bibitem{AC15a}
S.~N. Armstrong and P.~Cardaliaguet.
\newblock Quantitative stochastic homogenization of viscous {H}amilton-{J}acobi
  equations.
\newblock {\em Comm. Partial Differential Equations}, 40(3):540--600, 2015.

\bibitem{ACS14}
S.~N. Armstrong, P.~Cardaliaguet, and P.~E. Souganidis.
\newblock Error estimates and convergence rates for the stochastic
  homogenization of {H}amilton-{J}acobi equations.
\newblock {\em J. Amer. Math. Soc.}, 27(2):479--540, 2014.

\bibitem{AS131}
S.~N. Armstrong and P.~E. Souganidis.
\newblock Stochastic homogenization of {H}amilton-{J}acobi and degenerate
  {B}ellman equations in unbounded environments.
\newblock {\em J. Math. Pures Appl.}, 97(5):460--504, 2012.

\bibitem{AS132}
S.~N. Armstrong and P.~E. Souganidis.
\newblock Stochastic homogenization of level-set convex {H}amilton-{J}acobi
  equations.
\newblock {\em Int. Math. Res. Not. IMRN}, 2013(15):3420--3449, 2013.

\bibitem{AT141}
S.~N. Armstrong and H.~V. Tran.
\newblock {Stochastic homogenization of viscous Hamilton-Jacobi equations and
  applications}.
\newblock {\em Analysis \& PDE}, 7(8):1969--2007, 2014.

\bibitem{ATY15}
S.~N. Armstrong, H.~V. Tran, and Y.~Yu.
\newblock {Stochastic homogenization of a nonconvex Hamilton-Jacobi equation}.
\newblock {\em Calc. Var. Partial Differential Equations}, \noop{3003}in press.

\bibitem{ATY152}
Scott~N Armstrong, Hung~V Tran, and Yifeng Yu.
\newblock Stochastic homogenization of nonconvex hamilton--jacobi equations in
  one space dimension.
\newblock {\em Journal of Differential Equations}, 261(5):2702--2737, 2016.

\bibitem{BZZ06}
M.~Bramson, O.~Zeitouni, and M.~Zerner.
\newblock Shortest spanning trees and a counterexample for random walks in
  random environments.
\newblock {\em The Annals of Probability}, 34(3):821--856, 2006.

\bibitem{CIL92}
M.~G. Crandall, H.~Ishii, and P.-L. Lions.
\newblock User's guide to viscosity solutions of second order partial
  differential equations.
\newblock {\em Bull. Amer. Math. Soc. (N.S.)}, 27(1):1--67, 1992.

\bibitem{E92}
L.~C. Evans.
\newblock Periodic homogenisation of certain fully nonlinear partial
  differential equations.
\newblock {\em Proc. Roy. Soc. Edinburgh Sect. A}, 120(3-4):245--265, 1992.

\bibitem{F14}
B.~Fehrman.
\newblock A partial homogenization result for nonconvex viscous hamilton-jacobi
  equations.
\newblock {\em preprint arXiv:1402.5191}, 2014.

\bibitem{G15}
H.~Gao.
\newblock Random homogenization of coercive hamilton-jacobi equations in 1d.
\newblock {\em preprint arXiv:1507.07048}, 2015.

\bibitem{K07}
E.~Kosygina.
\newblock Homogenization of stochastic hamilton-jacobi equations: brief review
  of methods and applications.
\newblock {\em Contemporary Mathematics}, 429:189--204, 2007.

\bibitem{KRV06}
Elena Kosygina, Fraydoun Rezakhanlou, and Srinivasa~RS Varadhan.
\newblock Stochastic homogenization of hamilton-jacobi-bellman equations.
\newblock {\em Comm. Pure Appl. Math.}, 59(10):1489--1521, 2006.

\bibitem{LPV86}
P.-L. Lions, G.~Papanicolaou, and S.~R.~S. Varadhan.
\newblock Homogenization of hamilton-jacobi equations.
\newblock {\em Unpublished preprint}, 1986.

\bibitem{LS03}
P.-L. Lions and P.~E. Souganidis.
\newblock Correctors for the homogenization of hamilton-jacobi equations in the
  stationary ergodic setting.
\newblock {\em Comm. Pure Appl. Math.}, 56(10):1501--1524, 2003.

\bibitem{LS10}
P.-L. Lions and P.~E. Souganidis.
\newblock Stochastic homogenization of {H}amilton-{J}acobi and
  ``viscous''-{H}amilton-{J}acobi equations with convex
  nonlinearities---revisited.
\newblock {\em Commun. Math. Sci.}, 8(2):627--637, 2010.

\bibitem{LS05}
Pierre-Louis Lions and Panagiotis~E Souganidis.
\newblock Homogenization of “viscous” hamilton--jacobi equations in
  stationary ergodic media.
\newblock {\em Communications in Partial Difference Equations}, 30(3):335--375,
  2005.

\bibitem{MN12}
I.~Matic and J.~Nolen.
\newblock A sublinear variance bound for solutions of a random
  {H}amilton--{J}acobi equation.
\newblock {\em Journal of Statistical Physics}, 149(2):342--361, 2012.

\bibitem{RT00}
F.~Rezakhanlou and J.~E. Tarver.
\newblock Homogenization for stochastic {H}amilton-{J}acobi equations.
\newblock {\em Arch. Ration. Mech. Anal.}, 151(4):277--309, 2000.

\bibitem{S09}
R.~W. Schwab.
\newblock Stochastic homogenization of {H}amilton-{J}acobi equations in
  stationary ergodic spatio-temporal media.
\newblock {\em Indiana Univ. Math. J.}, 58(2):537--581, 2009.

\bibitem{S99}
P.~E. Souganidis.
\newblock Stochastic homogenization of {H}amilton-{J}acobi equations and some
  applications.
\newblock {\em Asymptot. Anal.}, 20(1):1--11, 1999.

\bibitem{Y89}
V.~Yurinskii.
\newblock On the error of averaging multidimensional diffusions.
\newblock {\em Theory of Probability \& Its Applications}, 33(1):11--21, 1989.

\end{thebibliography}

\end{document}